\documentclass[12pt]{article}

\usepackage{paralist}
\usepackage{natbib}
\setcitestyle{comma,aysep={}}
\usepackage{epsf}
\usepackage[pdftex]{graphicx}
\usepackage{amssymb,amsmath}
\usepackage{lscape}
\usepackage{umoline}
\usepackage{color}
\usepackage{multirow}
\usepackage{booktabs}
\usepackage{tikz}
\usetikzlibrary{shapes,arrows}
\usepackage[section]{placeins}
\usepackage[vlined]{algorithm2e}
\usepackage{secdot}

\usepackage{tikz}
\usetikzlibrary{positioning}
\usetikzlibrary{arrows}
\usetikzlibrary{shapes.arrows}
\usetikzlibrary{matrix}

\newcommand{\myfigure}[3]{ \begin{figure}[tbp]
    \begin{center}
    #1
    \end{center}
    \caption{#2} \label{#3} \vskip -0.14in \end{figure} }

\newcommand{\mytable}[3]{ \begin{table}[tbp]
  \begin{center}
  \caption{#2}
    #1
  \label{#3} \vskip -0.01in \end{center} \end{table} }

\begin{document}

\baselineskip 20pt plus .3pt minus .1pt

\begin{center}

{\LARGE Empirical analysis of metaheuristic search techniques for the parameterized dynamic slope scaling procedure }\\[12pt]

% Authors and addresses:

\footnotesize
\mbox{\large Weili Zhang, Charles D. Nicholson}\\
Department of Industrial and Systems Engineering, 
University of Oklahoma, 202 West Boyd, Room 449, Norman, OK, 73019-1022, USA 
\{\mbox{weili.zhang-1@ou.edu},\mbox{cnicholson@ou.edu} \}\\[6pt]

\normalsize

\end{center}

\noindent 

Abstract: The dynamic slope scaling procedure is an approximation method successfully which solves the fixed charge network flow (FCNF) problem by iteratively linearizing the fixed cost. The parameterized dynamic slope scaling procedure adds an additional $\psi$ parameter to the procedure which can significantly improve the solution quality. Finding the optimal value of $\psi$ for a given problem is non-trivial. This paper employs multiple metaheuristic techniques, including simulated annealing, tabu search, and particle swarm optimization, to guide the search for good parameter values.  In rigorous testing, we examine the search results, compare the improvement efficiencies among the techniques, and evaluate the final solution quality of the FCNF problem.  The experiments show that the solution improvement is robust with respect to these metaheuristics and the complexity of FCNF problem.

\bigskip

% Here are the Keywords:

% Here is the History:
\noindent {\it History:} TBD

\noindent\hrulefill
% The body of the paper starts here:

\section{Introduction}
\label{sec-Introduction}
The fixed charge network flow problem (FCNF) was first proposed by \citet{hirsch1954fixed}. The FCNF problem has many variants such as the  fixed charge transportation problem (FCTP) \citep{balinski1961fixed}, the lot sizing problem \citep{steinberg1980optimal}, and facility location problem \citep{nozick2001fixed}.  Each of these have common characteristics with the FCNF, namely the link between nodes in a given graph have both variable and fixed costs. The nodes have supply or demand amounts and the objective is to choose links and flow values on links to transfer commodities from supply nodes to demand nodes at a minimum cost. Given its wide spread application on many problem types \citep[e.g.][]{jarvis1978optimal,lederer1998airline,armacost2002composite, zhang2017bridge, zhang2018probabilistic, zhang2010lattice, zhang2016multi, zhang2010reformed, zhang2017resilience, zhang2016resilience}, the FCNF has been the focus of significant research for many years.

The FCNF problem is NP-hard and thus no efficient solution techniques exist which guarantee optimality \citep{GareyJ79}. Branch-and-bound (B\&B) \citep{land1960automatic} is commonly employed to solve the FCNF problem exactly \citep{driebeek1966algorithm, kennington1976new, barr1981new, cabot1984some, palekar1990branch}. B\&B is in general a very important algorithm in integer programing. The algorithm finds solutions by iteratively solving various relaxations of the problem to find integer feasible solutions.  At each iteration the lower bounds on the objective value (based on linear relaxations) and upper bounds (based on feasible solutions) may be updated.  As the gap between the bounds decreases, the B\&B approaches an optimal solution.  However, such a technique is inefficient in solving FCNF instances which lack tight bounds (e.g. due to high fixed costs).  Therefore, many researchers are interested in developing approximate methods for the FCNF problem. 

\citet{sun1998tabu} combined Tabu Search with the simplex on a graph method \citep{KennHelgason} to solve the FCTP. To help avoid local minima, their tabu implementation used both recency and frequency based memories and the simplex method was used to guide the .
\citet{adlakha2010heuristic} developed a heuristic algorithm based on Balinski's work \citeyearpar{balinski1961fixed}. This algorithm proposed a simple way to obtain valid upper bounds by solving related simple linear programs. This method is efficient with small sized problems.
\citet{monteiro2011ant} proposed a hybrid heuristic, using ant colony optimization to search broadly and as well as a local search method. 
\citet{molla2011solving} proposed an artificial immune algorithm and a genetic algorithm based on the spanning tree and pr\"{u}fer number representation to solve the FCTP in a two-stage supply chain network. 
\citet{antony2012genetic} applied genetic algorithm two scenarios of transportation problems, the first one was the FCTP and the second one considered opening costs of the distribution centers. Machine learning based approach is first developed by \citet{zhang2016prediction, nicholson2016optimal, zhang2018osea} which makes a significant contribution to the field.

\citet{KimPardalos99} developed the dynamic slope scaling procedure (DSSP) from the viewpoint of marginal cost and applied it solving the FCNF problem. The main idea is to linearize the objective function and to solve the problem iteratively, adapting the linearization at each step. DSSP removes the binary variables and incorporates the variable cost and fixed cost into the cost coefficient, which is updated by the previous linear programing solution. DSSP has been improved by combining with other methods : trust interval techniques \citep{kim2000dynamic}, intensification/diversification mechanisms \citep{crainic2004slope}, Tabu Search \citep{kim2006enhanced,gendron2003tabu} and parameterized DSSP \citep{NicholsonBarker2014parameterized}. Several types of FCNF variants, including concave piecewise linear network flow problem \citep{kim2000dynamic}, multicommodity location problem \citep{gendron2003tabu}, multicommodity FCNF \citep{crainic2004slope}, stochastic integer programming problem \citep{shiina2012dynamic}. 

This paper extends the work of \citet{NicholsonBarker2014parameterized} to explore a variety of metaheuristics search implementation for the parametrized dynamic slope scaling procedure ($\Psi$-DSSP).  The contribution in this paper is three-fold: (1) to explore the potential increase in the FCNF heuristic solution quality by using  more precise parameter values for $\Psi$-DSSP than have been examined to-date, (2) compare multiple well-known search metaheuristics applied to this problem with respect to solution quality and efficiency, and (3) evaluate solution efficacy by technique and parameter refinement across a spectrum of difficult FCNF instances.
The remainder of this paper is organized as follows: Section \ref{sec-Background} briefly provides the background knowledge of fixed charge network flow problem and parameterized dynamic slope scaling procedure; Section \ref{sec-MetaSearch} details three metaheuristics techniques and the corresponding implementations with respect to $\Psi$-DSSP; Section \ref{sec-Experiment} describes the experimental design to verify the new approach; Section \ref{sec-Analysis} analyzes the results from experiments and Section \ref{sec-Conclusions} summarizes the findings.

\section{Background}
\label{sec-Background}
\subsection{FCNF Problem Description}
The single-commodity, uncapacitated, fixed charge network flow problem can be defined on a directed graph $\mathcal{D}=(N,A)$, where $N$ is the set of nodes and $A$ is the set of arcs $(i,j)$. Each node $i \in N$ has a commodity requirement value $R_i$ ($R_i>0$ for supply nodes, $R_i<0$ for demand nodes, $R_i=0$ for transshipment nodes). Each arc $(i,j)\in A$ has an associated fixed and variable cost $f_{ij}$ and $c_{ij}$, respectively.  The decision variable $x_{ij}$ denotes the flow on arc $(i,j)\in A$. The fixed cost on arc $(i,j)$ is the cost of using the arc for any positive flow, $x_{ij} >0$.  The binary decision variables $y_{ij}$ for all $(i,j) \in A$ are used to model which arcs are selected for commodity flow.  The aim of the FCNF is to select a subset of arcs to be opened and determine the flow on the arcs such that the supply in the network is routed to meet the demand at a minimal total cost. The FCNF is formulated as a traditional mixed binary programing problem as follows:
\begin{align} 
\text{min} & \sum_{(i,j) \in A} (c_{ij}x_{ij} + f_{ij}y_{ij})  \label{eq1} \\
\text{s.t.} & \sum_{(i,j) \in A} x_{ij} - \sum_{(j,i) \in A} x_{ji} = R_{i}&
  \forall{i \in N} \label{eq2} \\
&0 \leq x_{ij} \leq M_{ij}y_{ij}&
  \forall{(i,j) \in A} \label{eq3} \\
&y_{ij} \in \lbrace 0,1 \rbrace&
  \forall{(i,j) \in A} \label{eq4}
\end{align}

\noindent  The objective in \eqref{eq1} is nonlinear. The constraints in \eqref{eq2} ensure that the inflow and outflow of a node satisfy the supply/demand at node $i \in N$. 
%$c_{ij}$ denotes the variable cost of arc  $(i,j)\in A$;
%$f_{ij}$ denotes the fixed cost of arc $(i,j)\in A$;
%$R_i$ denotes the supply/demand at node $i \in N$;
$M_{ij}$ is an arc capacities (artificial or real) are used in the constraints \eqref{eq3} to ensure that the fixed cost $f_{ij}$ is incurred whenever there is a positive flow on arc $(i,j)$.  If any arc $(i,j)$ is not capacitated, then $M_{ij}$ can be set to a sufficiently large value to not inhibit feasible flow.  The constraints in \eqref{eq4} define $y_{ij}$ as binary for all $(i,j) \in A$.

\subsection{ Parameterized Dynamic Slope Scaling Procedure}
The dynamic slope scaling procedure replaces the variable cost and fixed cost with a  linearized cost coefficient, which is updated iteratively by the marginal fixed cost. Let $\tilde{x}^{k}_{ij}$ denote the flow value of arc $(i,j)$ at solution iteration $k$,
\begin{align*}
\tilde{x}^k_{ij} = 
\begin{cases}
x^k_{ij} & x^k_{ij} > 0 \\
\tilde{x}^{k-1}_{ij} & x^k_{ij} = 0
\end{cases}
\qquad \forall(i,j) \in A.
\end{align*}

\noindent Let $\bar{c}_{ij}^{k+1}$ denote the linearized cost coefficient of arc $(i,j)$ at iteration $k+1$, 
\begin{align*}
\bar{c}^{k+1}_{ij}=c_{ij}+\frac{f_{ij}}{\tilde{x}^{k}_{ij}}
\end{align*}
The FCNF objective can be formulated at each iteration $k$:
\begin{align*} 
\text{min} &  \sum_{(i,j) \in A} \left(c_{ij} + \frac{f_{ij}}{\tilde{x}^{k-1}_{ij}}\right)x^k_{ij}  
\end{align*}
The algorithm stops when no more improvements can occur, or equivalently when
\begin{align*}
\tilde{x}^{k-1}_{ij} = \tilde{x}^k_{ij} \qquad\forall(i,j) \in A.
\end{align*}
\noindent The objective value associated with the original problem formulation \eqref{eq1} is computed at each iteration.  The objective associated with the best solution across all DSSP iterations is returned (e.g. not necessarily the last iteration).

\citet{NicholsonBarker2014parameterized} provided solid evidence that a small modification of the linearizion technique can considerably alter the algorithm's search path and significantly enhance the solution quality over the DSSP. The revised algorithm is parameterized by a value $\psi$ such that the corresponding problem formulation for the $k^{\text{th}}$ iteration of $\Psi$-DSSP is
\begin{align}
\hspace{1em} \text{min} &  \sum_{(i,j) \in A} \left(c_{ij} + \psi\frac{f_{ij}}{\tilde{x}^{k-1}_{ij}}\right)x^k_{ij}  \\
\text{s.t.} & \sum_{(i,j) \in A} x^{k}_{ij} - \sum_{(j,i) \in A} x^{k}_{ji} = R_{i} &   \forall{i \in N}  \\
            & 0 \leq x^{k}_{ij} \leq M_{ij} &  \forall{(i,j) \in A} 
\end{align}
In their work they examined various values of $\psi\in(0,2)$ and ultimately tested solution quality on difficult FCNF problems for five fixed values, $\psi = \{0.25, 0.5, 0.75, 1, 1.25\}$. The results were promising with maximum improvements of $\Psi$-DSSP over DSSP nearing 30\%.  The paper establishes the sensitivity of solution quality with respect to the $\psi$ parameter and also demonstrates there is no clear pattern relating $\psi$ with solution efficacy for different problem instances. 

\section{Metaheuristics for Searching $\psi$ Values}
\label{sec-MetaSearch}
Since the best $\psi$ values vary among instances and affect the solution  significantly, in this study metaheuristics are employed as a search procedure to find good $\psi$ values which improve the solution quality of the FCNF problem in reasonable time.  While several techniques are possible, we choose two classic single-solution based metaheuristic algorithms, simulated annealing \citep{SA1983} and tabu search \citep{glover1986future}, and  one population-based metaheuristic algorithm, particle swarm optimization \citep{kennedy1995particle}. The details of each implementation are discussed below.  For notational convenience, let $z = $ DSSP($\psi$) denote the best objective value from all iterations in the $\Psi$-DSSP approach for the parameter value equal to $\psi$.

\subsection{Simulated Annealing}
Simulated annealing (SA) is a well known and effective stochastic search algorithm useful on non-linear and discontinuous problems spaces \citep{SA1983,SA2003}.  While SA was originally developed for discrete problems, it has been successfully extended to continuous optimization problems \citep[e.g.][]{GenSA1986,ContinuousSA91}.  This technique is inspired by the physical annealing process in which a crystalline solid is heated and cooled slowly in such a way to improve the structural integrity of the material.  Essentially, at the beginning of the process, the high temperatures allows the atoms to move freely and move to a state of minimum energy during the cooling process. Similarly for optimization problems, SA allows for a diverse search at the beginning of the process which may include very poor solutions.  This diversity allows the process to escape local minima.  At later iterations, as the ``temperature'' decreases, the search becomes more restrictive and focused on solutions which improve the objective value.  For this investigation we will employ the \emph{Boltzmann annealing} \citep{van1987simulated} and \emph{very fast annealing} \citep{Ingber1989} implementations to find appropriate values for $\psi$. A brief explanation of the procedure follows.                   

Let $\psi_{0}$ denote the initial value for the DSSP parameter and $z_{0}$ denote the corresponding best solution value for that parameter setting.  Let the value $T_{0}$ denote the initial temperature in the simulated annealing method.  At each iteration $i$, choose a value $\psi_{c} \in \left[\psi_{i}-\frac{\epsilon}{2}, \psi_{i}+\frac{\epsilon}{2} \right]$ as a candidate parameter in the ``neighborhood'' of $\psi_{i}$, where $\epsilon >0$ is the range of the neighborhood.  (The method for choosing $\psi_{c}$ from the neighborhood is detailed below.)  Evaluate the candidate parameter and compare the solution quality with that of the solution found using $\psi_{i}$.  If $z_{c} \le z_{i}$, then the candidate parameter is ``accepted'' by setting   $\psi_{i+1} = \psi_{c}$.  Otherwise, the candidate solution is accepted with probability 
\begin{align}
p=\exp \left[ \frac{   1 - \left( \frac{z_{c}}{z_{i}} \right)  }{T_{i}} \right]  \label{SAaccept}          
\end{align}
\noindent where $T_{i}$ is updated by different functions based on the SA annealing implementation:
\begin{align}
\text{Boltzmann: }  &  T_{i} =\frac{T_{0}}{\log(1+i)} \label{tempScheduleB}\\
\text{Very Fast Annealing: } &  T_{i} = T_{0}e^{\left(\frac{-i}{e}\right)}.
\label{tempScheduleVFA}
\end{align}
\noindent In equation \eqref{SAaccept}, the probability computed is proportional to the percent change between $z_{c}$ and $z_{i}$.  This ensures the acceptance criteria for a given temperature setting is the same for a candidate objective value which is 10\% higher than current objective, regardless of the magnitude of the objectives.  Also note in equations \eqref{tempScheduleB} and \eqref{tempScheduleVFA} various factors are available to further tune the search procedure.  We have set these factors to 1 for simplicity.  

If the candidate solution is not accepted, the procedure evaluates up to a pre-specified number of other candidate values (denoted $i_{dwell}$) before continuing to the next iteration and decreasing the temperature value. Selecting a value $\psi_{c}$ from the neighborhood of $\psi_{i}$ is executed  differently depending on the implementation.  For $u_{i}$ chosen randomly and uniformly on $(0,1)$ and for $v_i$ randomly selected from the standard normal distribution, the neighbor values are chosen as follows:
\begin{align*}
\text{Boltzmann: } & \psi_{c} = \begin{cases}
\psi_{i} + \frac{1}{2}v_i\sqrt{T_{i}} \qquad & \text{if } \sqrt{T_{i}} < \frac{2}{3}\\[2ex]
\psi_{i} + \frac{1}{3}v_i \qquad & \text{otherwise}
\end{cases}\\
\text{Very Fast Annealing: } & \psi_{c} = \begin{cases}
\psi_{i} + T_{i}\left[\left(1+\frac{1}{T_i}\right)^{|2u_i-1|}-1\right] \qquad & \text{if } u_i < 0.5\\[2ex]
\psi_{i} - T_{i}\left[\left(1+\frac{1}{T_i}\right)^{|2u_i-1|}-1\right] \qquad & \text{otherwise}
\end{cases}\\
\end{align*}
%\noindent The temperature updating scheme and the method for choosing candidate solutions is based on \citet{Ingber1989}.

\noindent Pseudo code for $\Psi$-DSSP combined with simulated annealing is shown in Figure \ref{SADSSP}.  The stopping criterion are based on maximum iterations ($i_{max}$) or maximum total runtime ($t_{max}$).

\myfigure{
\begin{algorithm}[H]
\DontPrintSemicolon
%\KwData{$\mathcal{D}=(N,A), \mathbf{f,c,M,R}, i_{max}, t_{max},  T_0, \psi_0$}
\KwData{$\Psi$-DSSP instance, $i_{max}, t_{max},  T_0, \psi_0$}
\KwResult{$\psi_{best},z_{best}$}
\Begin{
 % $z_0 \gets$ objective value of $\Psi$-DSSP for $\psi_0$\;
 $z_0 \gets$ DSSP($\psi_0$)\;
 $ z_{best} \gets z$, $\psi_{best} \gets \psi_0$ \;
  $i \gets 0, t \gets 0$\;
  \While{$i < i_{\text{max}}$  {\bf and}  $t < t_{max}$ {\bf and} early termination \emph{is false}} {
    $i \gets i + 1$\;
    $j \gets 0, t \gets$ total runtime\;
    update $T_i$\;
    
    \While{$j < i_{dwell}$  {\bf and}  $\psi_i \neq \psi_c$}{
    $j \gets j + 1$\;
    $\psi_{c} \gets $ choose neighbor of $\psi_{i}$\;
   % $z_{c} \gets$ objective value of $\Psi$-DSSP for $\psi_c$\; 
     $z_c \gets$ DSSP($\psi_c$)\;
     \eIf{ $z_{c} < z_{best} $ }{
      $z_{best} \gets z_{c}$,$\psi_{best} \gets \psi_{c}$ \;
      $ z_{i} \gets z_{c}$, $\psi_{i} \gets \psi_{c}$ \;
      }{
      generate $u \sim U(0,1)$\;
      $p \gets \exp \left[ \frac{   1 - \left( \frac{z_{c}}{z_{i}} \right)  }{T_{i}} \right]$ \;
      \If{$u<p$}{
      $\psi_{i} \gets \psi_{c}$ }
     }
   }
  Set \emph{early termination} boolean value
  }
  \Return{$\psi_{best}, z_{best}$} 
  }
\end{algorithm}
}{Simulated Annealing for $\Psi$-DSSP Implementation}{SADSSP}

\subsection{Tabu Search}
Tabu Search (TS) was developed by \citet{glover1977heuristics,glover1986future, glover1989tabu,glover1990tabu}.  The algorithm improves on  basic \emph{hill-climbing} (a greedy method which always seeks improved solutions) by allowing the search to move away from good solutions in attempt to escape local optima. Once the search moves to a new solution in the neighborhood, the previous solution is added to a \emph{tabu list} which forbids the search to return the previous states for a certain number of iterations.  TS is an effective approach for combinatorial problems such as graph coloring \citep{dubois1993epcot} and the linear ordering problem \citep {DuarteLaguna2011}. The algorithm has been extended to continuous optimization problems by developing a paradigm of concentric hyperspheres \citep{siarry1997fitting} or hypercubes \citep{chelouah2000tabu} to define the neighborhood of a solution. 

The implementation of TS for the $\Psi$-DSSP problem is for a one dimensional, continuous search space.  We employ the neighborhood definition of \citet{siarry1997fitting} as a series of $k$ concentric intervals.  For each iteration $i$, the neighborhood of $\psi_{i}$ are the $k$ concentric bands $B = \{B_1,B_2,\ldots,B_k\}$, where    
\begin{align}
B_j = \{\psi_{c}|h_{j-1}\leq \parallel \psi_{c}-\psi_{i}\parallel \leq h_{j}\} 
\label{Interval}
\end{align}
and $h_0,h_1,\dots,h_k$ is an increasing positive sequence of values. The value for $h_k$ is chosen as the maximum absolute difference between a given solution and its furthest neighbor.  The values $h_j$ are determined according to a geometric series
\begin{align*}
h_{k-j}=\frac{h_k}{2^{j}} \qquad \text{for } j=1,2,\ldots, k-1
\end{align*}

\myfigure{
\begin{tikzpicture}[scale=1]

%begin graph on left ---------
\footnotesize
% line
\draw[thick] (-6.5,0) -- (6.5,0); % node[anchor=north] {$M_{ij}$};

% ticks and labels
\draw[-] (-6,-0.1) -- (-6,0.1);
\draw	(-6,-0.1) node[anchor=north] {$-h_3$};

\draw[-] (-3,-0.1) -- (-3,0.1) ; %node[anchor=east] {$z_{ij}$};
\draw	(-3,-0.1) node[anchor=north] {$-h_2$};

\draw[-] (-1.5,-0.1) -- (-1.5,0.1) ; %node[anchor=east] {$z_{ij}$};
\draw	(-1.5,-0.1) node[anchor=north] {$-h_1$};

\draw[-] (-0.25,-0.1) -- (-0.25,0.1) ; %node[anchor=east] {$z_{ij}$};
\draw	(-.3,-0.1) node[anchor=north] {$-h_0$};

\draw[-] (6,-0.1) -- (6,0.1) ; %node[anchor=east] {$z_{ij}$};
\draw	(6,-0.1) node[anchor=north] {$h_3$};

\draw[-] (3,-0.1) -- (3,0.1) ; %node[anchor=east] {$z_{ij}$};
\draw	(3,-0.1) node[anchor=north] {$h_2$};

\draw[-] (1.5,-0.1) -- (1.5,0.1) ; %node[anchor=east] {$z_{ij}$};
\draw	(1.5,-0.1) node[anchor=north] {$h_1$};

\draw[-] (0.25,-0.1) -- (0.25,0.1) ; %node[anchor=east] {$z_{ij}$};
\draw	(.3,-0.1) node[anchor=north] {$h_0$};

%add neighbor points
\draw [fill] (0,0) circle [radius=0.075];
\draw (0,0.1) node[anchor=south] {$\psi_i$};

\draw [fill] (1,0) circle [radius=0.075];
\draw (1,0.1) node[anchor=south] {$\psi_{i,1}$};

\draw [fill] (-1.85,0) circle [radius=0.075];
\draw (-1.85,0.1) node[anchor=south] {$\psi_{i,2}$};

\draw [fill] (5.5,0) circle [radius=0.075];
\draw (5.5,0.1) node[anchor=south] {$\psi_{i,3}$};

\end{tikzpicture}

}
{Partition of the $\psi_{i}$ neighborhood when $k=3$}{interval} 

 %The external radii, $h_{k}$, and internal radii,

\noindent The value $h_0$ is the smallest positive value in the sequence, but is not related directly to $h_k$. It is used in the algorithm to specify the minimum distance between two $\psi$ values for those two values to be considered distinct.  

In iteration $i$, the $k$ neighbors of $\psi_i$ are selected by randomly choosing one point from each $B_j$ for $j=\{1,2,\ldots,k\} $.  An example with $k=3$ is depicted in Figure \ref{interval}.  The neighbor value $\psi_{i,1}$ is selected from the set of intervals closest to $\psi_i$, the value $\psi_{i,2}$ is selected from the second closest set of intervals, and so on. Pseudo code for $\Psi$-DSSP combined with tabu search is provided in Figure \ref{TSDSSP}.

\myfigure{
\begin{algorithm}[H]
\DontPrintSemicolon
\KwData{$\Psi$-DSSP instance, $i_{max}, t_{max},k,n, h_0,h_k, \psi_0$}
\KwResult{$\psi_{best},z_{best}$}
\Begin{
   $z_0 \gets$ DSSP($\psi_0$)\;
   
   $z_{best} \gets z_0$, $\psi_{best} \gets \psi_0$ \;
   
   $i \gets 0, t \gets 0$\;
   
   $tabu list \gets$ Null \;
  
   \While{$i < i_{\text{max}}$  {\bf and}  $t < t_{max}$ {\bf and} early termination \emph{is false}}
   {
      Form concentric bands $\{B_1,\ldots,B_k\}$ around $\psi_i$ using $h_0,\ldots,h_k$ \;
      \For {$j=1,\ldots,k$}
      { 
         randomly select $\psi_{i,j}$ from $B_j$ s.t. $\parallel \psi_{i,j} - \psi_{tabu} \parallel > h_0$ $\forall  \psi_{tabu} \in tabu list$\;
      }
    $index \gets \underset{j \in \{1,2,\ldots,k\}}{\arg\min}$ DSSP($\psi_{i,j}$) \;
    
    $z_{bestcandidate} \gets$ DSSP($\psi_{i,index}$)\;
   
    $\psi_{i} \gets \psi_{i,index}$ \;
   
   \If{$z_{bestcandidate} < z_{best}$}
   {
     $z_{best} \gets z_{bestcandidate}$ \;
     $\psi_{best} \gets \psi_{i}$ \;
   }
   add $\psi_{i}$ to end of $tabu list$ \;
   \If{ $\left\vert{tabu list}\right\vert > n$ }{
       remove the first tabu element }
    $i \gets i+1, t \gets$ total runtime\; 	
    Set \emph{early termination} boolean value
    }
  \Return{$\psi_{best}, z_{best}$} 
  }
\end{algorithm}
}{Tabu Search for $\Psi$-DSSP Implementation}{TSDSSP}

\subsection{Particle Swarm Optimization}
Artificial life \citep{adami1998introduction} simulates natural biotic system behavior with the help of computer or other abiotic system \citep{wilke2002biology}. One main branch of artificial intelligence is evolutionary computing \citep{back1997handbook} inspired by Darwin's theory of evolution. Another important branch is swarm intelligence \citep{bonabeau1999swarm} inspired by social group behavior. Swarm optimization is similar to evolutionary algorithms (e.g. genetic algorithm) in that it is a stochastic population-based metaheuristic. Swarm optimization differs in that the same population persists throughout many iterations and the members of the population adapt their behavior based on their own history and from ``learning'' from other members. In particle swarm optimization (PSO)  \citep{kennedy1995particle,eberhart1995new}, particles``fly'' through the solution space analogous to how birds flock or fish swarm. In PSO, each particle has a position vector, velocity vector, and fitness value. The particle's position is a point in the solution space, the velocity determines (with stochastically) which solution the particle will move to next, and the fitness vis the objective value of the current solution. Each particle maintains a record of the value and position of its individual best historical fitness value. The also swarm maintains the global best historical fitness value and position. The algorithm uses these pieces of information to inform all particles and update their positions.  Figure \ref{pso} graphically depicts the particle influences during a single iteration.  The final updated position derives from a linear combination of influences from its own iteration history, the swarm history, and the particle's current velocity vector.  PSO has been empirically shown to outperform genetic algorithms with respect ot speed of convergence and in some cases quality of the solutions \citep{angeline1998evolutionary,hassan2005comparison}. Furthermore, PSO has fewer parameters to refine than competing evolutionary algorithms which makes this approach even more appealing.

%\myfigure{
%\includegraphics[scale=0.5]{pso.png}
%}{Particle Swarm Optimization}{pso}

\myfigure{
\begin{tikzpicture}[scale=1]

%begin graph on left ---------
\footnotesize
% line
\draw[thick] (-6.5,0) -- (6.5,0); 

\draw[dotted] (-1.5,-1) -- (-1.5,0.75); 

%add psi points points
\draw [fill] (-1.5,0) circle [radius=0.1];
\draw (-2.0,0.1) node[anchor=south] {$\psi_{u,i}$};

\draw [fill] (0.75,0) circle [radius=0.1];
\draw (0.95,0.1) node[anchor=south] {$\psi_{u,i+1}$};

\draw [fill] (3.55,0) circle [radius=0.1];
\draw (3.55,0.1) node[anchor=south] {$\psi^{\text{pbest}_u}$};

\draw [fill] (-5.0,0) circle [radius=0.1];
\draw (-5.0,0.1) node[anchor=south]  {$\psi^{\text{gbest}}$};

%add influence lines  
\draw[dashed,-triangle 45] (-1.5,-.5) -- (3.55,-.5);                    %particle best
\draw (1.25,-1.1) node[anchor=south] {particle best influence};

\draw[-triangle 45,line width=0.75mm] (-1.5,0) -- (0.75,0);                %final

\draw[dashed,-triangle 45] (-1.5,.75) -- (2,.75);
\draw (-.25,.75) node[anchor=south] {$v_{u,i}$};

\draw[dashed,-triangle 45] (-1.5,-1) -- (-5,-1);                    %swarm best
\draw (-2.75,-1.5) node[anchor=south] {swarm influence};

\end{tikzpicture}

}
{Particle Swarm Optimization on a Line}{pso} 

% If $z_{u,i} < z^{pbest_u}$, then the best parameter value for particle $u$ is updated: $z^{\text{pbest}_u}=z_{u,i}$ and $\psi^{\text{pbest}_u} = \psi_{u,i}$. If the $z_{u,i} < z^{\text{gbest}}$, where  $z^{\text{gbest}}$ is the swarm best found fitness value, then the swarm best value and swarm best parameter are updated by setting $ z^{\text{gbest}}=z_{u,i}$ and $\psi^g = \psi_{u,i}$. 

Let $n$ denote the total number of particles and $\psi_{u,i}$ denote the position of particle $u$ in iteration $i$. Let $i_{max}$ denote the maximum number of iterations. In iteration $i$, the fitness of particle $u$ is $z_{u,i}$ and is computed by solving  $\Psi$-DSSP with $\psi_{u,i}$.  The best historical position and objective for a given particle is maintained as $\psi^{\text{pbest}_u}$ and $z^{\text{pbest}_u}$, respectively.  Similarly,  $\psi^{\text{gbest}}$ and $z^{\text{gbest}}$ denote the global best historical values for particle positions and objective values for the swarm.  Let $v_{u,i}$ denote the velocity of particle $u$ at iteration $i$ and at the end of iteration $i$, $v_{u,i}$ is updated according to
\begin{align}
w_{i} & =w_{max}-\frac{i}{i_{max}} \left(w_{max}-w_{min} \right)
\label{updatew}\\
v_{u,i+1} & = w_{i}v_{u,i}+c_1r_1(\psi^{\text{pbest}_u}-\psi_{u,i})+c_2r_2(\psi^{\text{gbest}}-\psi_{u,i})
\label{updateposition}
\end{align}
where $w_{max}, w_{min}, c_1$, and $c_2$ are control the ``inertia'' and ``acceleration'' of the particles and are useful in tuning the search. The inertia parameter $w_{i}$ in \eqref{updatew} decreases linearly  with respect to the total number of permissible iterations, beginning at $w_{max}$ and ending at $w_{min}$.  This encourages a diverse search during the beginning iterations and convergence in the later iterations.  The values for $c_1$ and $c_2$ are constants which reflect the movement influence associated with a given particle's historical best position and the swarm's best historical position. The parameters $r_1$ and $r_2$ are two uniformly random numbers between $0$ and $1$. Then the updated position for particle $u$ is computed by equation(\ref{updatepsi}),
\begin{align}
\psi_{u,i+1}=\psi_{u,i}+v_{u,i+1}
\label{updatepsi}.
\end{align}

\noindent Pseudo code for $\Psi$-DSSP combined with particle swarm optimization to solve FCNF is provided in Figure \ref{PSODSSP}.
\myfigure{
\begin{algorithm}[H]
\DontPrintSemicolon
\KwData{$\Psi$-DSSP instance, $i_{max}, t_{max}, n, w_{max},w_{min},v_{max},c_1,c_2$}
\KwResult{$\psi^{\text{gbest}},z^{\text{gbest}}$}
\Begin{
	\For{$u=1,\ldots,n$}{
	    generate $\psi_{u,0} \sim U(0,2)$ \;
	    generate $v_{u,0} \sim U(-v_{max},v_{max})$ \;
	    $\psi^{\text{pbest}_u} \gets \psi_{u,0}$\;
	    $z^{\text{pbest}_u} \gets \infty$\;
	    }
   $z^{\text{gbest}} \gets \infty$,$\psi^{\text{gbest}} \gets$ Null \;  
   $i \gets 0, t \gets 0$\;
   \While{$i < i_{max}$  {\bf and}  $t < t_{max}$  {\bf and} early termination \emph{is false}}
   {
   	\For{$u=1,\ldots,n$}
   	{
   	$z_{u,i} \gets \text{DSSP}(\psi_{u,i})$\;
   	\If{$z_{u,i} \leq z^{\text{pbest}_{u}}$}
   	{
   	 $z^{\text{pbest}_u} \gets z_{u,i}$ \;
   	 $\psi^{\text{pbest}_u} \gets \psi_{u,i}$\;
   	 \If{$z_{u,i} \leq z^{\text{gbest}}$}
   	 { 
   	 $z^{\text{gbest}} \gets z_{u,i}$ \;
   	 $\psi^{\text{gbest}} \gets \psi_{u,i}$\;
   	 }
   	}
   	} 
   	$w_{i} \gets w_{max}-\frac{i(w_{max}-w_{min})}{i_{max}}$\;
   	\For{$u=1,\ldots,n$}
   	{
   	generate $r_1,r_2 \sim U(0,1)$\;
   	$v_{u,i+1} \gets w_{i}v_{u,i}+c_1r_1(\psi^{\text{pbest}_u}-\psi_{u,i})+c_2r_2(\psi^{\text{gbest}}-\psi_{u,i})$\;
   	$\psi_{u,i+1} \gets \psi_{u,i}+v_{u,i+1}$\;
   	}
   $i \gets i+1, t \gets$ total runtime\; 	
   Set \emph{early termination} boolean value\;
  \Return{$\psi^{g}, z^{g}$} 
  }
  }
\end{algorithm}
}{Particle Swarm Optimization for $\Psi$-DSSP Implementation}{PSODSSP}

\section{Computational Experiments Design}
\label{sec-Experiment}
\subsection{Network Characteristics}
In order to evaluate the efficiency and solution quality from the various metaheuristic approaches to the $\Psi$-DSSP problem, the experimental design includes tests on a variety of network sizes, each which has characteristics corresponding to difficult FCNF problem instances (e.g. high fixed to variable cost ratio).  The tests includes networks with 25, 50, and 100 nodes.  For each level of node quantity $n$, we generate 60 feasible FCNF instances in which the number of arcs $m$, is randomly selected. Specifically, we randomly choose $n-1 \leq \frac{m}{2} \leq \frac{n(n-1)}{2}$ and create a connected network instance where each of the $\frac{m}{2}$ undirected arcs is replaced by two directed arcs.

%To provide a solid foundation for metaheuristics searching for good $\psi$ values, we design the test with a larger variety of nodes ($25$,$50$,$100$) and solve $60$ instances within each set of nodes. The number of arcs $m$ is randomly generated between a connected graph and full connected.

According to \citet{NicholsonBarker2014parameterized}, the gap between the na\"{i}ve DSSP objective and optimal value is related to the network characteristics.  In our experimentation we focus on the instances in which DSSP performed poorly as identified in their work.  We also use the same network characteristics measures and instance specification as \citet{NicholsonBarker2014parameterized}. The percentage of supply, demand, and transshipment nodes are respectively randomly selected with approximate probabilities $0.2$, $0.2$ and $0.6$.  The probabilities are approximate in that adjustments are made to ensure an instance is feasible.  The variable costs and fixed costs for each link are randomly assigned on $U(0,20)$ and $U(20000,60000)$, respectively. The total requirements for each supply node is randomly assigned on $U(1000,2000)$. The total requirements of supply node is distributed randomly as negative requirements to the demand nodes. Gurobi 5.6 is used to solve the linear programing problems. 

To describe the network characteristics let $d=\frac{m}{2\binom{n}{2}}$ denote the density of the network, $S$ denote the total supply, and $\rho_s$ and $\rho_d$ represent the percentage of supply and demand nodes, respectively.  Let $\theta=\frac{S}{n_s}$ indicate the average supply for each supply node. Let $\phi$ denote the overall network ratio of fixed to  variable costs, 
\begin{align*}
\phi &= \frac{\displaystyle\sum_{(i,j) \in A} f_{ij}}{\displaystyle\sum_{(i,j) \in A} c_{ij}}
\end{align*}
where $f_{ij}$ and $c_{ij}$ denote the fixed cost and variable cost of arc $(i,j)$. Let  $\gamma=\frac{\theta}{\phi}$ be a characteristic designed to provide an apriori estimate of the ratio of the cost components of a feasible solution. That is, since $\theta$ is proportional to the magnitude of arc flow in a feasible solution, then
\begin{align*}
\gamma &= \frac{\theta}{\phi} = \frac{ \displaystyle  \sum_{(i,j) \in A} \theta  c_{ij}}{\displaystyle\sum_{(i,j) \in A} f_{ij}}
\end{align*}
\noindent is defined so that if $\gamma\gg1$, variable costs are likely more important than fixed costs; and if $\gamma\ll1$, then fixed costs are likely to have a larger influence on the optimal solution. 

\mytable{
\begin{center}
\begin{tabular}[!ht]{c c c c c c}
\toprule
Variable	 & 	mean	 & 	std.dev	 & 	min	 & 	max	 & 	skewness	\\
\midrule
$d$	 & 	0.5	 & 	0.3	 & 	0.03	 & 	1.0	 & 	0.2	\\
$\rho_s$	 & 	31.8\%	 & 	6.4\%	 & 	20.0\%	 & 	56.0\%	 & 	0.4	\\
$\rho_d$	 & 	20.3\%	 & 	5.2\%	 & 	4.0\%	 & 	36.0\%	 & 	0.2	\\
$\theta$	 & 	1020.9	 & 	342.0	 & 	118.8	 & 	1760.2	 & 	0.2	\\
$\phi$	 & 	3994.7	 & 	136.8	 & 	3552.1	 & 	5093.7	 & 	2.7	\\
$\gamma$	 & 	0.3	 & 	0.1	 & 	0.0	 & 	0.5	 & 	0.2	\\
\bottomrule
\end{tabular}
\end{center}}{Network Characteristics} {characteristics}

Table \ref{characteristics} summarizes the network characteristics (mean, standard deviation, minimum, maximum, and skewness) for the test instances.  Figure \ref{histogramSPhiRhoGamma} shows the corresponding histograms of each characteristic.  The least dense instance ($d=0.03$) is composed of 100 nodes and 316 arcs, whereas the most dense network is composed of 100 nodes and 9690 arcs.  The test cases include network instances with  6\% to 56\% of nodes as supply nodes, and from 5\% to 36\% as demand nodes.  The supply quantity per supply node ranges from about 120 to nearly 1,800.  The fixed to variable network cost ratio averages to about 4,000.  The value for $\gamma$ ranges from nearly 0 to 0.5. 

\myfigure{
\centering
\includegraphics[scale=.8]{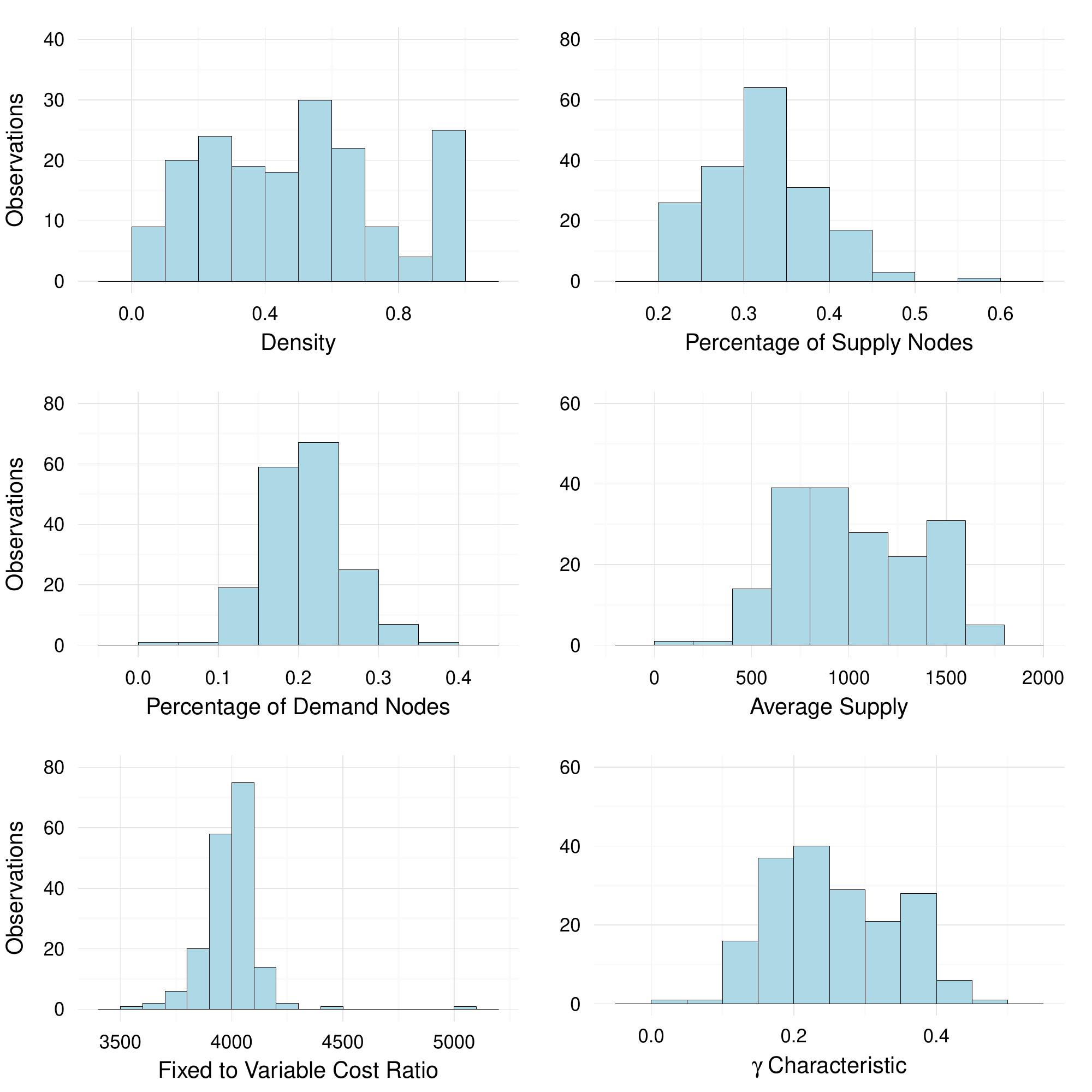}
}{Distribution of Network Characteristics}{histogramSPhiRhoGamma}

\subsection{Parameters Setting in Metaheuristics}
For any metaheuristic implementation, the specific tuning parameters and stopping criterion play a significant role in the efficacy of the technique.  Our testing does not evaluate the full range of the manifold possibilities, however we have attempted to define reasonable and comparable settings for the techniques we use.  Three stopping criterion are employed in all the techniques.  We classify these as hard criteria or early termination criteria.  The hard criterion is twofold: (i) maximum time limit $t_{max}$ (1 hour) and (ii) the maximum number of allowable iterations $i_{max}$ (which is set equal to $t_{max}$ divided by the time it takes to solve the network instance using na\"{i}ve DSSP). The $i_{max}$ and $t_{max}$ are two hard stop criteria, but we can not compare the efficiency of different algorithms by these criteria since the best solution maybe found before meeting $i_{max}$ or $t_{max}$.  An early termination criterion is employed which stops the search when the objective value is not improved for several iterations.  After extensive experimentation and suggestions from literature in section \ref{sec-MetaSearch}, we set the metaheuristic parameters as listed in Table \ref{settings}.

\mytable{
\begin{tabular}{l l}
Metaheuristic & Parameter Settings \\
\toprule
\multirow{2}{*}{Simulated Annealing} & max iterations at a given temperature: $i_{dwell} = 3$ \\
& initial temperature: $T_0 = 0.25$ \\ 
\midrule
\multirow{4}{*}{Tabu Search} & minimum neighbor distance: $h_0 = 0.01$ \\
                             & maximum neighbor distance: $h_k = 0.2$ \\
                             & maximum length of $tabulist$: $n=5$ \\
                             & quantity of neighbors:  $k=5$\\
\midrule
Particle Swarm Optimization & number of particles:  $n= 10$\\
                            & maximum velocity: $v_{max}=1$ \\
							& inertia weight: $w_{min}=0.4,w_{max}=0.9$\\
							& acceleration coefficients: $c_1=c_2=2$\\
\bottomrule
\end{tabular}}{Metaheuristic Settings}{settings}

\section{Results Analysis}
\label{sec-Analysis}
In this section, we report the correlation between solution quality performance and network characteristics, correlation between $\psi$ value and improvement, and the time efficiency of the different approaches: simulated annealing using the Boltzmann annealing (SAB), simulated annealing using ``very fast annealing'' (SAVF), tabu search on continuous space using concentric bands (TS), and particle swarm optimization (PSO). Figure \ref{psi} shows the histograms of the best empirical $\psi$ values across all $180$ instances for each method.  The best identified parameter value spans nearly the entire search space on $[0,2]$.
  
  \myfigure{
  \includegraphics[scale=.8]{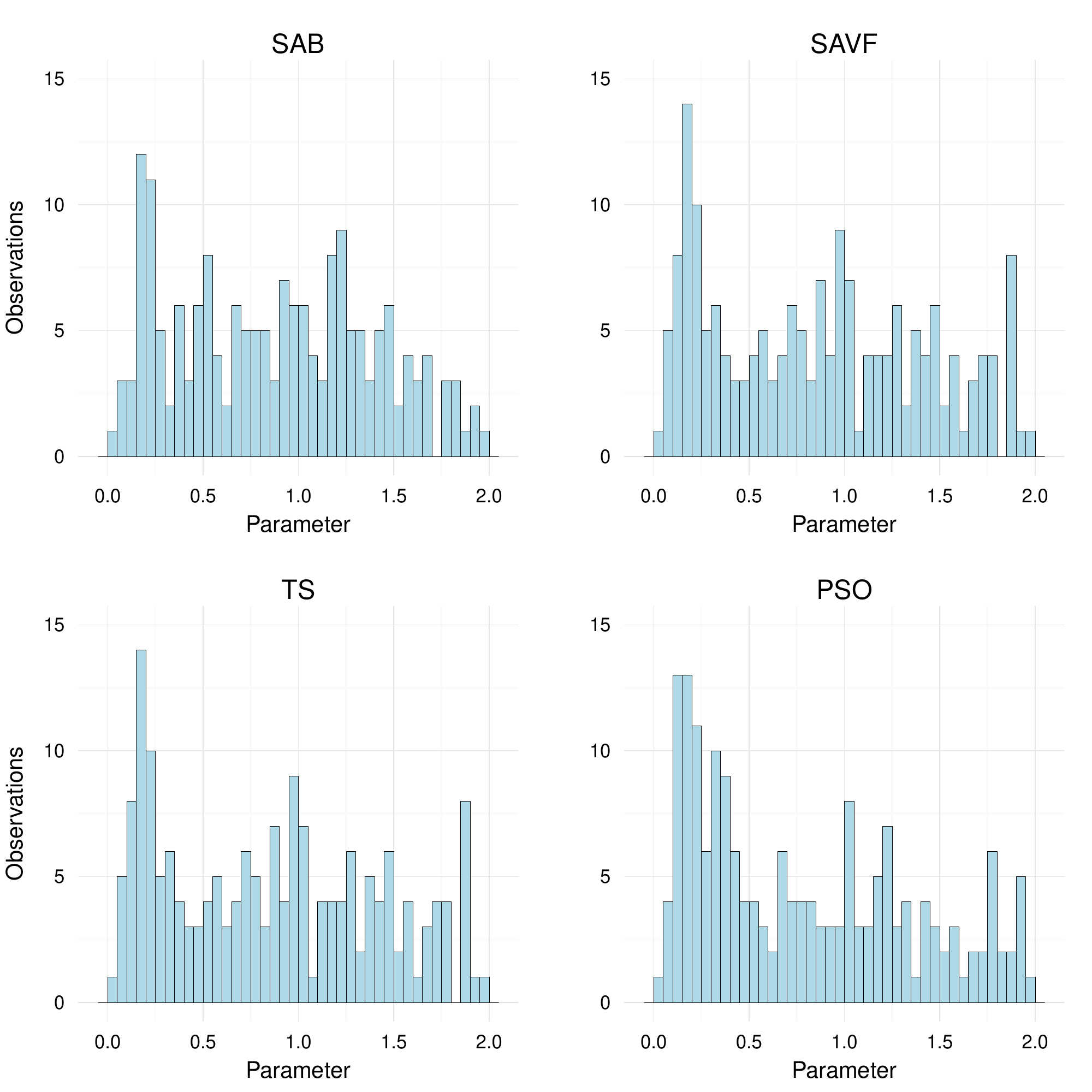}
  }{Distribution of Best Parameter Value}{psi}
  
Let $z_{\text{DSSP}}$ denote the best objective value resulting from the na\"{i}ve DSSP, and $z_{\text{SAB}}$,  $z_{\text{SAVF}}$,  $z_{\text{TS}}$,  $z_{\text{PSO}}$ denote the best objective value resulting from $\Psi$-DSSP and the associated search technique.  Let $z^{\text{gap}}$ denote the percentage improvement of a given test approach compared to DSSP.  That is, $z^{\text{gap}}$ is the percent decrease in the objective value
\begin{align*} 
z^{gap} &= \frac{z_{\text{DSSP}} - z_{\text{x}}}{z_{\text{DSSP}}} \times 100\%.
\end{align*} 
\noindent where $z_{\text{x}}$ is one of $z_{\text{SAB}}$,  $z_{\text{SAVF}}$,  $z_{\text{TS}}$, or $z_{\text{PSO}}$. The parameterized DSSP using metaheuristic search outperforms na\"{i}ve DSSP in 680 (94.37\%) cases with an overall mean improvement of 11.83\%.  Figure \ref{histgap} shows the gap distribution for $z^{\text{gap}}>0$.  \color{red} The improvement percentage ranges from $0.2\%$ to $24.2\%$.   24\% is the maximum.... \color{black}

\myfigure{
\includegraphics[scale=1]{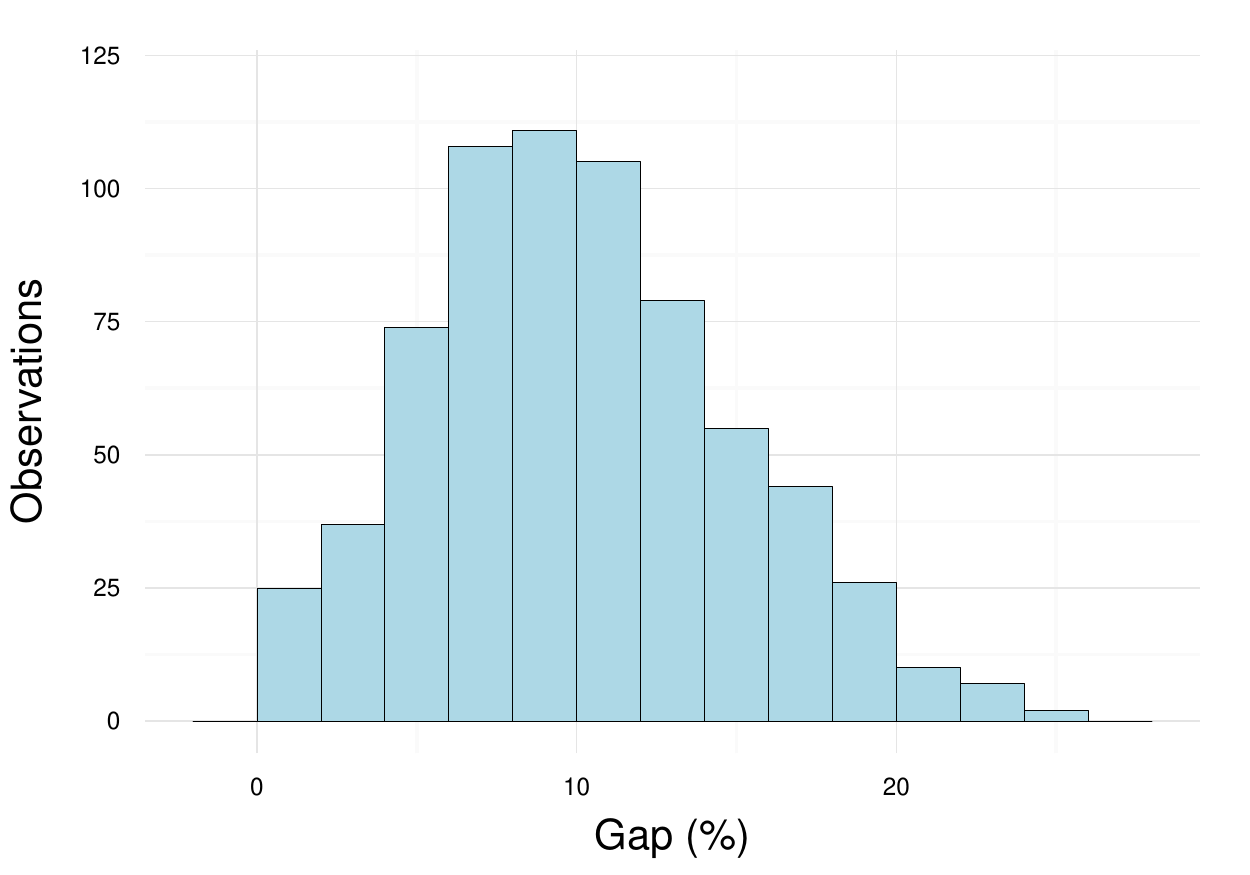}
}{Distribution of Gap}{histgap}

\myfigure{
\centering
\includegraphics[scale=1]{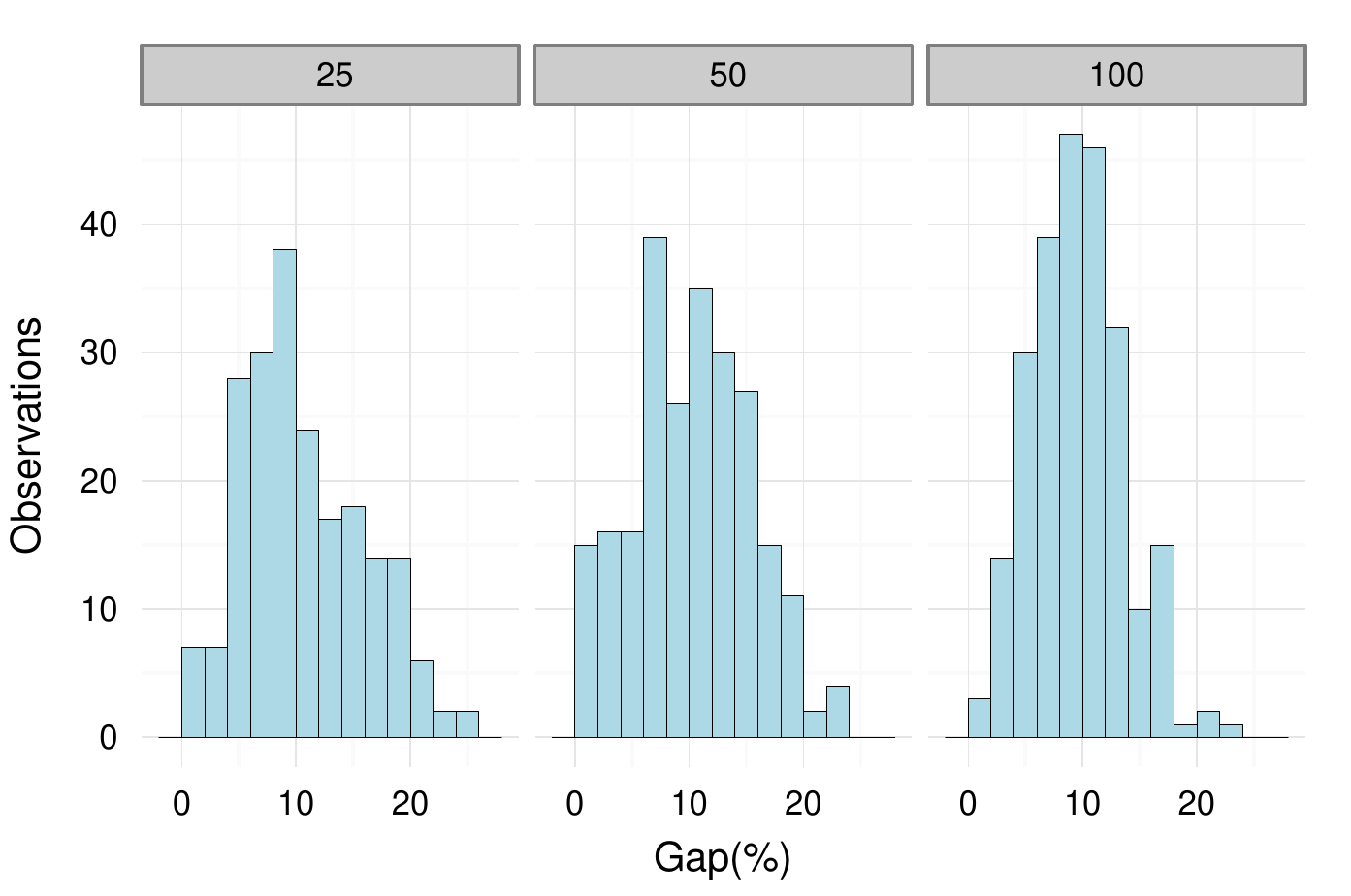}
}{Distribution of Gap by Number of Nodes}{histgapbynode}

Table \ref{gapstatistics} reports the percentage of instances in which $\Psi$-DSSP outperforms DSSP, and the corresponding mean, min, and max $z^{\text{gap}}$ overall and for each metaheursitic search approach for each level of nodes for the FCNF instances.  The results between different metaheuristics are quite close. The average improvement overall and for each individual approach is around 10\%. The percentage of cases in which $\Psi$-DSSP with metaheuristics outperforms DSSP increases with the number of nodes and achieves $100\%$ when the number of nodes equals to $100$. However, the magnitude of the gap itself for instances with $z^{\text{gap}} > 0$ shows no statistical difference with respect to the number of nodes ($p$-value = 0.91).  Furthermore, there is no significant difference in the average solution gap by search technique.  Figure \ref{histgapbynode} shows the gap distribution by number of nodes for the problem instances with improvement. %From , we can see that when number of nodes equal to $100$, the distribution is more centric to $10\%$. Therefore, the $\Psi$-DSSP with metaheuristics is recommended to use when the problem is larger and more complex. We can also identify that the gap is not influenced by the number of nodes, which can be verified by the ANOVA test (p-value $=0.91$).

%The minimum gap of $100$ nodes is $10$ times larger than the minimum gap of $25$ nodes. 

\mytable{
\begin{tabular}{c c c c c c c}
\toprule
 & & \multicolumn{3}{c}{Nodes} & \\ \cmidrule{3-5}
	& & 	25	 & 	50	 & 	100	 & 	Overall	\\
\midrule
\multirow{3}{*}{Overall} & $z^{\text{gap}} > 0$ (\%)	 & 	86.0	 & 	98.3	 & 	100.0	 & 	94.9	\\
& mean $z^{\text{gap}}$ 	 & 	10.7	 & 	10.3	 & 	9.6	 & 	10.2	\\
& max $z^{\text{gap}}$ 	 & 	24.2	 & 23.9	 & 	22.0	 & 	24.2	\\
\midrule
\multirow{3}{*}{SAB} & $z^{\text{gap}}_\text{SAB} > 0$ (\%)	 & 	86.7	 & 	98.3	 & 	100.0	 & 	95.0	\\
& mean $z^{\text{gap}}_\text{SAB}$ 	 & 	10.7	 & 	10.5	 & 	9.6	 & 	10.2	\\
& max $z^{\text{gap}}$  	 & 23.9	 & 	23.9	 & 	17.0	 & 	23.9	\\
\midrule
\multirow{3}{*}{SAVF} & $z^{\text{gap}}_\text{SAVF} > 0$ (\%)	 & 	86.7	 & 	98.3	 & 	100.0	 & 	95.0	\\
& mean $z^{\text{gap}}_\text{SAVF}$ 	 & 	10.4	 & 	10.0	 & 	10.0	 & 	10.1	\\
& max $z^{\text{gap}}$  	 & 	24.2	 & 	22.5	 & 	22.0	 & 	24.2	\\
\midrule
\multirow{3}{*}{TS} & $z^{\text{gap}}_\text{TS} > 0$ (\%)	 & 	86.7	 & 	98.3	 & 	100.0	 & 	95.0	\\
& mean $z^{\text{gap}}_\text{TS}$ 	 & 	10.8	 & 	10.1	 & 	9.4	 & 	10.1	\\
& max $z^{\text{gap}}$  	 & 	24.2	 & 	23.7	 & 	20.3	 & 	24.2	\\
\midrule
\multirow{3}{*}{PSO} & $z^{\text{gap}}_\text{PSO} > 0$ (\%)	 & 	85.0	 & 	98.3	 & 	100.0	 & 	94.4	\\
& mean $z^{\text{gap}}_\text{PSO}$ 	 & 	10.9	 & 	10.5	 & 	9.5	 & 	10.3	\\
& max $z^{\text{gap}}$  	 & 	23.1	 & 	23.9	 & 20.3	 & 	23.9	\\
\bottomrule
\end{tabular}
}{Gap Statistics} {gapstatistics}

\mytable{
\begin{center}
\begin{tabular}[!ht]{ccccccccc}
  \toprule
    \multirow{2}{*}{Variable} & \multicolumn{2}{c}{$z^{\text{gap}}_{\text{SAB}}$}  & \multicolumn{2}{c}{$z^{gap}_\text{SAVF}$} & \multicolumn{2}{c}{$z^{gap}_\text{TS}$} & \multicolumn{2}{c}{$z^{gap}_\text{PSO}$}\\
    \cmidrule(lr){2-3} \cmidrule(lr){4-5} \cmidrule(lr){6-7} \cmidrule(lr){8-9} & Cor & $p$-value & Cor   & $p$-value & Cor   & $p$-value & Cor   & $p$-value \\
  \midrule
  $d$	 & 	0.31	 & 	1.86E-05	 & 	0.26	 & 	4.87E-04	 & 	0.23	 & 	1.80E-03	 & 	0.24	 & 	1.07E-03	\\
  $\rho_s$	 & 	-0.09	 & 	0.21	 & 	-0.10	 & 	0.20	 & 	-0.07	 & 	0.33	 & 	-0.10	 & 	0.19	\\
  $\rho_d$	 & 	-0.06	 & 	0.40	 & 	-0.02	 & 	0.78	 & 	-0.10	 & 	0.20	 & 	-0.09	 & 	0.24	\\
  $\theta$	 & 	0.01	 & 	0.94	 & 	-0.05	 & 	0.48	 & 	-0.03	 & 	0.68	 & 	-0.04	 & 	0.59	\\
  $\phi$	 & 	-0.05	 & 	0.54	 & 	-0.08	 & 	0.30	 & 	-0.06	 & 	0.39	 & 	-0.07	 & 	0.37	\\
  $\gamma$	 & 	0.01	 & 	0.90	 & 	-0.05	 & 	0.54	 & 	-0.03	 & 	0.73	 & 	-0.04	 & 	0.64	\\
  \bottomrule
\end{tabular}
\end{center}}{Correlation Analysis} {cor}
The Pearson correlations listed in Table \ref{cor} show that only density of the network has a statistically significant positive correlation with $z^{gap}$. As the density increases, the complexity of the FCNF rises and the $z^{gap}$ increases regardless of the search method employed. No significant correlations exist between the gap and the remaining descriptors.  That is, with the exception of network density, the improvement attainable from using search techniques is relatively robust with respect to the range of network characteristics evaluated  in this study.  There is also no correlation between $z^{\text{gap}}$ and the parameter value itself (pearson correlation coefficient $=0.004$). Figure \ref{psigap} plots the 720 $\psi$ values and solution gap pairs.  The best $\psi$ value varies with different problem instances.

\myfigure{
\includegraphics[scale=1]{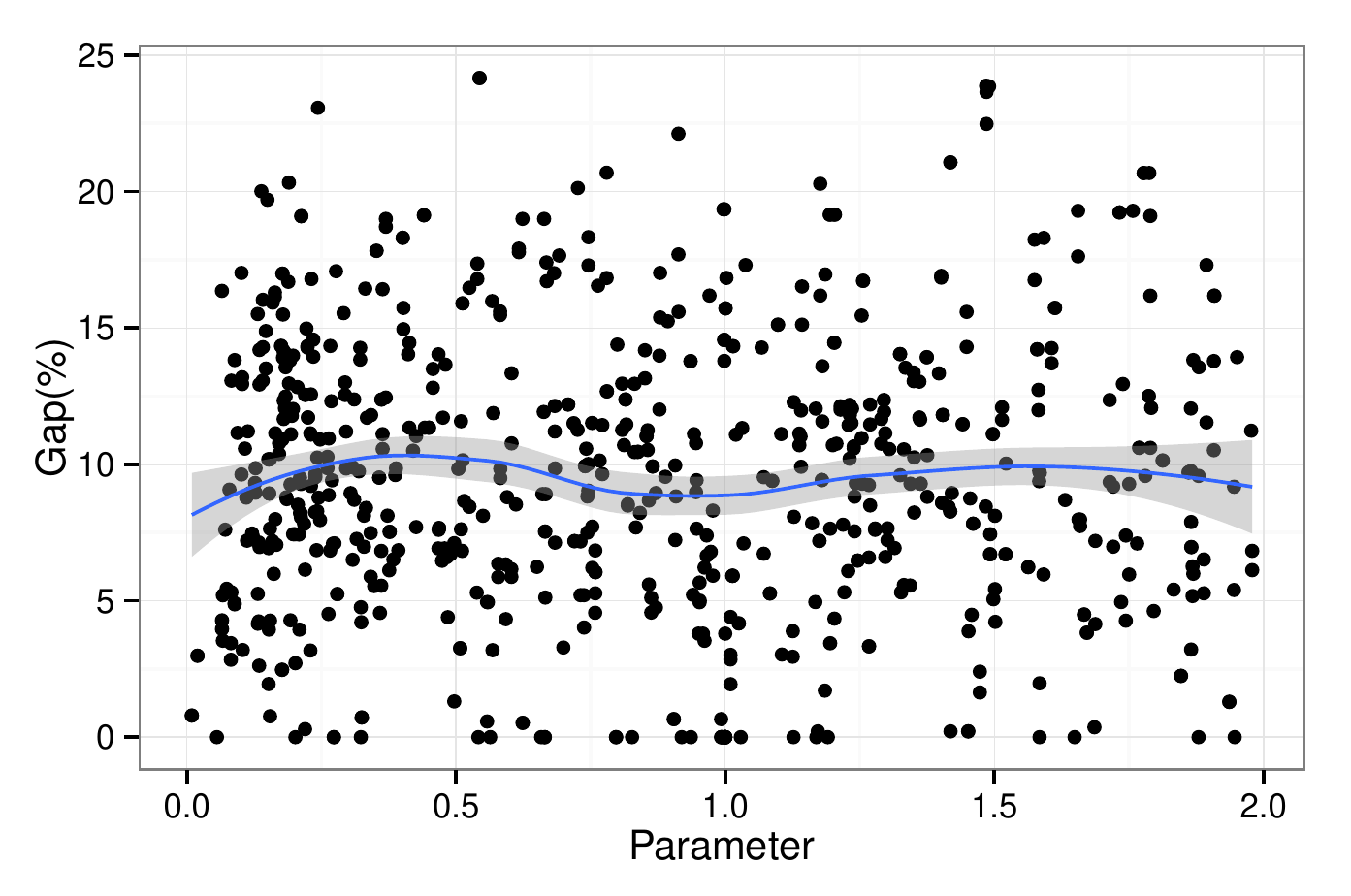}
}{Scatter Plot of $\psi$ Value and Gap}{psigap}

The various search implementations with $\Psi$-DSSP each produce notable improvements over DSSP.  To compare the techniques we consider three elements: $z_\text{gap}$, CPU time, and total iterations. Each technique was allowed to run for a maxmimum of 1 hour, however some criteria allow for early termination.  Terminating early is a good quality assuming the solution is also good.  However, early termination due to quick convergence to relatively poor local optimum is not a desirable feature.  To account for this we all report the average solution quality per algorithm iteration.  Specifically, let $s$ denote the total number of parameterized DSSP solutions evaluated for a given search procedure.  Let $r_x$ represent the solution efficiency defined as the solution gap per iteration, $r_x = \frac{z^{gap}_x}{s_x}$, and $t_x$ denote the CPU time for each search technique $x \in \{ \text{SAB, SAVF, TS, PSO}\}$. For example, a typical $\Psi$-DSSP with PSO search might use 10 particles with 100 updates each, resulting in 1,000 DSSP solutions. If this procedure produced a gap of 15\%, then $r_\text{PSO} = \frac{15}{1000} = 0.015$.

Table \ref{anova} reports the average CPU time and average gap per iteration ($r_x$) for each of the solution approaches for each level of nodes.  The Tukey HSD tests results are reported if the $p$-value of the ANOVA test is smaller than $0.05$.  For the problem instances with 25 nodes, the swarm optimization outperformed the other methods, although the average time for all methods was relatively small.   Statistically, all four methods had identical performance for the problem instances with 50 nodes.  For the largest problem sizes, on average both PSO and TS converged after about 30 minutes of runtime.  This was significantly faster than either of the two simulated annealing implementations (48 minutes for SAB; 1 hour and 9 minutes for SAVF).

%\mytable{
%\begin{tabular}{cccc}
%  \toprule
%    \multirow{2}{*}{Node Level}	 & 	 \multirow{2}{*}{Metric}	 &  Mean & 	 \multirow{2}{*}{Tukey Result} \\
%                                 &                               & SAB,SAVF,TS,PSO &  
%     \\
%  \midrule
%    \multirow{3}{*}{$25$} & $z^{gap}$	& 11,10,11,11 & 	no significant difference	\\
%  & 	$t$	& 14,16,11,7  & 	$t_\text{PSO} < t_\text{TS} < t_\text{SAB} = t_\text{SAVF}$
% \\
%   & $r$($\times 10^{-4}$)	 & 22,17,25,43	 & 	$r_\text{PSO}>r_\text{TS}>r_\text{SAB}=r_\text{SAVF}$
%  \\
%  \midrule
%      \multirow{3}{*}{$50$}  & 	$z^{gap}$	& 11,10,10,11	 & 	no significant difference	\\
%       & $t$	& 249,282,166,249	 & 	no significant difference	\\
%       &  $r$($\times 10^{-4}$)	& 8,7,9,8 	 & 	no significant difference	\\
%      
%  \midrule
%      \multirow{3}{*}{$100$}  & 	$z^{gap}$	& 10,10,9,10	 & 	no significant difference	\\
%       & 	$t$	 & 2889,4114,1934,1597  & 	$t_\text{PSO}=t_{TS}<t_\text{SAB}<t_\text{SAVF}$	\\
%       &  $r$($\times 10^{-4}$)	 & 3,2,3,5	 & 	$r_\text{PSO}>r_\text{TS}=r_\text{SAB}=r_\text{SAVF}$	\\
%      
%  \bottomrule
%\end{tabular}
%}{ANOVA and Tukey Test} {anova}

\mytable{
\begin{tabular}{ccccccc}
  \toprule
    \multirow{2}{*}{Node Level}	 & 	 \multirow{2}{*}{Metric}	 &  \multicolumn{4}{c}{Mean} & 	 \multirow{2}{*}{Tukey Result} \\ \cmidrule{3-6}
                                 &                               & SAB & SAVF & TS & PSO &  
     \\
  \midrule
    \multirow{2}{*}{$25$} & 	CPU time (sec)	& 14 & 16& 11 & 7  & 	$t_\text{PSO} < t_\text{TS} < t_\text{SAB} = t_\text{SAVF}$
 \\
   & $r$ ($\times 10^{-4}$)	 & 22 & 17 & 25 & 43	 & 	$r_\text{PSO}>r_\text{TS}>r_\text{SAB}=r_\text{SAVF}$
  \\
  \midrule
      \multirow{2}{*}{$50$}   & 	CPU time (sec)	& 249 &  282 & 166 & 249	 & 	no significant difference	\\
       &  $r$ ($\times 10^{-4}$)	& 8 & 7 & 9 & 8 	 & 	no significant difference	\\
      
  \midrule
      \multirow{2}{*}{$100$}  & 	 	CPU time (sec) & 2889 & 4114 & 1934 & 1597  & 	$t_\text{PSO}=t_{TS}<t_\text{SAB}<t_\text{SAVF}$	\\
       &  $r$ ($\times 10^{-4}$)	 & 3 & 2 & 3 & 5	 & 	$r_\text{PSO}>r_\text{TS}=r_\text{SAB}=r_\text{SAVF}$	\\
      
  \bottomrule
\end{tabular}
}{ANOVA and Tukey Test} {anova}

%Second, considering the efficiency, PSO is the best choice in all the instances. When the number of nodes equals to $25$, PSO is significantly faster than the rest methods. When the number of nodes equals to $50$, these methods could be replaced equivalently. When the number of nodes equals to $100$, even though the $t^{PSO}$ is equivalent to $t^{TS}$, the $s^{PSO}$ is significant larger than $s^{TS}$. The running time determined by the iterations and computation capacity of computers, the iteration is an absolute metrics to measure the efficiency of metaheuristics algorithms. Hence, PSO is the best choice among all the tested techniques. 

\section{Conclusions}
\label{sec-Conclusions}
The FCNF problem is a classically NP-hard problem with many real-world applications.  Due to it's complexity, many techniques have been developed to approximate solutions quickly.  The dynamic slope scaling procedure and the more general parameterized variation, $\Psi$-DSSP are two such approaches. In this study we address three questions with respect to these two techniques.  First, using a highly refined parameter setting, we quantify reasonable expectations for solution improvement of $\Psi$-DSSP over DSSP.  Empirically we find solution quality to be improved by a statistically significant 10\% on average, and up to 24\% as a maximum across various sized network problems.  Approximately 95\% of the instances benefit from a refined parameter value. 

Secondly, we evaluate the performance of multiple metaheuristic search algorithms with regard to  this problem class. Implementations of simulated annealing, tabu search, and particle swarm optimization all ultimately find equally good parameter values for $\Psi$-DSSP in our testing. In terms of efficiency, there are differences.  PSO on average converges faster than the competing techniques (e.g. about half of the time as simulated annealing) and has the highest average ``improvement per iteration''.  

Finally, we find that the solution quality of $\Psi$-DSSP is relatively robust with respect to a wide spectrum of FCNF network characteristics for each search technique.  Each method is tested on 180 different network instances with various quantities of nodes, arcs, commodity supplies and demands, and values for fixed and variable costs.  Arc density is the only feature with a  statistically significant correlation to solution quality, i.e. as the density increases, there is evidence that $\Psi$-DSSP further outperforms na\"{i}ve DSSP.  

We conclude that employing a metaheuristic strategy, especially PSO, as a search technique to accompany $\Psi$-DSSP is a reasonable approach.  There is no obvious relationship between network characteristics and an optimal $\psi$ value for the parameterized DSSP.  However, the improvement in solution quality is notable and most complex FCNF instances benefit from a refined parameter setting.  
    
\newpage

\bibliographystyle{ormsv080}
\bibliography{references}

\end{document}